\documentclass[12pt]{article}
\usepackage{amsmath,amsthm,amsfonts,latexsym,amssymb,enumerate,color}
\usepackage{graphicx}

\def\AA{\mathbb A}
\def\C{\mathcal C}
\def\CC{\mathbb C}
\def\DD{\mathbb D}
\def\E{\mathcal E}
\def\G{\mathcal G}
\def\M{\mathcal M}
\def\N{\mathcal N}
\def\NN{\mathbb N}
\def\RR{\mathbb R}
\def\TT{\mathbb T}
\def\ZZ{\mathbb Z}

\def\H2p{\overline{H^2_0}}

\def\beq{\begin{equation}}
\def\eeq{\end{equation}}

\def\BP{Blaschke product }

\def\clos{\mathop{\rm clos}\nolimits}

\def\ker{\mathop{\rm ker}\nolimits}
\def\Kmin{\mathop{\rm K_{min}}\nolimits}
\def\LCM{\mathop{\rm LCM}\nolimits}

\def\spam{\mathop{\rm span}\nolimits} 

\newtheorem{thm}{Theorem}[section]

\newtheorem{cor}[thm]{Corollary}

\newtheorem{ex}[thm]{Example}

\numberwithin{equation}{section}

\def\beginpf{\begin{proof}}
\def\endpf{\end{proof}}

\begin{document}

\title{Toeplitz kernels and model spaces}


\author{M.~Cristina C\^amara\thanks{
Center for Mathematical Analysis, Geometry and Dynamical Systems,
Instituto Superior T\'ecnico, Universidade de Lisboa, 
Av. Rovisco Pais, 1049-001 Lisboa, Portugal.
 \tt ccamara@math.ist.utl.pt} \ and  Jonathan R.~Partington\thanks{School of Mathematics, University of Leeds, Leeds LS2~9JT, U.K. {\tt j.r.partington@leeds.ac.uk}}\ \thanks{Corresponding author}
}

\maketitle

\begin{abstract}
We review some classical and more recent results concerning kernels of Toeplitz operators and their relations with model spaces, which are themselves Toeplitz kernels of a special kind. We highlight the fundamental role played by the existence of maximal vectors for every nontrivial Toeplitz kernel.
\end{abstract}

\noindent {\bf Keywords:}
Toeplitz kernel,  model space, nearly-invariant subspace,\\ minimal kernel, multiplier, Carleson measure

\noindent{\bf MSC:} 47B35, 30H10.

\maketitle

\section{Introduction}
We shall mostly be discussing Toeplitz operators on the Hardy space
$H^2=H^2(\DD)$ of the unit disc $\DD$, which embeds isometrically as a closed subspace
of $L^2(\TT)$, where $\TT$ is the unit circle, by means of non-tangential limits.
These are standard facts that can be found in many places, such as
\cite{duren,nik}.

For a symbol $g \in L^\infty(\TT)$ the Toeplitz operator $T_g :H^2 \to H^2$ is defined
by
\[
T_g f = P_{H^2} (g \cdot f) \qquad (f \in H^2),
\]
where $P_{H^2}$ denotes the orthogonal projection from $L^2(\TT)$ onto $H^2$.

Similarly we may define Toeplitz operators on the Hardy space $H^2(\CC^+)$ of the
upper half-plane, which embeds as a closed subspace of $L^2(\RR)$, and we shall
use the same notation, since the context should always be clear, writing
\[
T_g f = P_{H^2(\CC^+)} (g \cdot f ) \qquad (f \in H^2(\CC^+)),
\]
where $P_{H^2(\CC^+)}$ is the orthogonal projection from $L^2(\RR)$ onto $H^2(\CC^+)$.

The kernels of such operators have been a subject of serious study for at least fifty years, and one
particular example here is the class of model spaces.
Let $\theta \in H^\infty=H^\infty(\mathbb D)$ be an inner function, that is $|\theta(t)|=1$ almost everywhere on $\mathbb T$,
and consider the Toeplitz operator $T_{\overline\theta}$. It is easily verified that its kernel
is the space
\[
K_\theta:= H^2 \ominus \theta H^2 = H^2 \cap \theta \overline{H^2_0},
\]
where $\overline{H^2_0}$ denotes the orthogonal complement of $H^2$ in $L^2(\RR)$.
It follows from Beurling's theorem that these spaces $K_\theta$ 
are the nontrivial closed invariant subspaces of the backward shift operator $S^*=T_{\bar z}$, defined by
\[
S^* f(z) = \frac{f(z)-f(0)}{z} \qquad (f \in H^2, \quad z \in \DD).
\]
They
include the spaces of polynomials of degree at most $n$ 
for $n=0,1,2,\ldots$ (take
$\theta(z)=z^{n+1}$), as well as the finite-dimensional spaces consisting of rational
functions (each such $n$-dimensional space corresponds to taking $\theta$ to be a Blaschke product of degree $n$).
For a good recent book on model spaces, see  \cite{GMR}.

Another example, which has applications in systems and control theory, is the space
corresponding to the inner function $\theta_T(s)=e^{isT}$ in $H^\infty(\CC^+)$, for a fixed $T>0$.
For by the Paley--Wiener theorem, the Fourier transform establishes a canonical isometric isomorphism between
$L^2(0,\infty)$ and $H^2(\CC^+)$, mapping the subspace $L^2(0,T)$ onto $K_{\theta_T}$.

As we shall now see, the class of Toeplitz kernels, which includes the class of model spaces, can itself be
described in terms of model spaces. Most of the results we present are valid (with suitable modifications)
in $H^p$ for $1<p<\infty$, as well as in Hardy spaces on the half-plane. The interested reader may
refer back to the original sources.

We recall first one classical result of Coburn \cite{coburn}, that for $g \in L^\infty(\TT)$ not almost everywhere $0$,
either $\ker T_g=\{0\}$ or $\ker T^*_g=\{0\}$ (note that $T^*_g=T_{\overline g}$).
This was proved as an intermediate step towards showing that the Weyl spectrum of
a Toeplitz operator coincides with its essential spectrum.

\section{Background results}

\subsection{The 1980s}

The papers of Nakazi \cite{nak}, Hayashi \cite{hay85,hayashi86},  Hitt \cite{hitt} and Sarason \cite{sarason88}
were all published within a short space of time. 

Nakazi's paper is mostly concerned with finite-dimensional Toeplitz kernels, but does explore the role of 
{\em rigid\/} functions in the context of Toeplitz kernels.
He uses the term {\em $p$-strong} for an outer function $f \in H^p$ with the property that if $kf \in H^p$
for some measurable $k$ with $k \ge 0$ a.e., then $k$ is constant, although nowadays the term {\em rigid\/}
is generally adopted.
He then shows that $\dim\ker T_g = n$, a non-zero integer, if and only if $\ker T_g=u P_{n-1}$,
where $u \in H^2$ with $u^2$   rigid, and 
$P_{n-1}$ is the space of polynomials of degree at most $n-1$.
Nakazi's work also bears on extremal problems and the properties of Hankel operators.

In fact, a function $f \in H^1$ with $\|f\|=1$ is rigid if and only if it is an {\em exposed point} of the ball of $H^1$;
that is, if and only if there is a functional $\phi \in (H^1)^*$ such that
\[
\phi(f)=\|\phi\|=\|f\|=1,
\]
and such that if $\phi(g)=1$ for some $g$ with $\|g\|=1$, then $g=f$. Chapter 6 of \cite{FricMash}
contains a useful discussion of this result.\\

Meanwhile, Hayashi \cite{hayashi86} showed that
the kernel of a Toeplitz operator $T_g$ can be
written as $uK_{\theta}$, where $u$ is outer and $\theta$ is inner with $\theta(0)=0$, and $u$ multiplies the model space
$K_{\theta}$
isometrically  onto $\ker T_g$. 
Every closed subspace $M$ of $H^2$ possesses a reproducing kernel $k_w \in M$ (where $w \in \DD$), such that
$\langle f,k_w \rangle=f(w)$ for $f \in M$, and, as
an application of his main result, Hayashi gave an expression for the reproducing kernel corresponding to
a Toeplitz kernel, namely,
\[
k_w (z) = \overline{u(w)}u(z) \frac{1-\overline{ \theta(w)}\theta(z)}{1-\overline{w} z},
\]
for $w,z \in \DD$, where $\ker T_g = uK_{\theta}$.
Hayashi also noted in \cite{hay85} that every nontrivial  Toeplitz kernel $T_g$ is equal to
$\ker T_{\overline h/h}$ for some outer function $h$, a significant simplification in the analysis of
Toeplitz kernels. Moreover, in the representation $uK_\theta$, we have that $u^2$ is rigid.\\

Hitt's work was mostly concerned with the Hardy space $H^2(\AA)$ of the annulus $\AA=\{z \in \CC: 1 < |z| < R \}$
for some $R>1$, and in
classifying those
closed subspaces of $H^2(\AA)$ invariant under $Sf(z)=zf(z)$.
To do this he made a study of subspaces $M$ of $H^2(\DD)$ that are 
{\em nearly invariant\/}
 under
the backwards shift $S^*$, i.e., $f \in M$ and $f(0)=0$ implies that $S^*f \in M$. (Again, his original terminology,
{\em weakly invariant}, has been superseded.) 

It is easy to see that a Toeplitz kernel is nearly $S^*$-invariant, for if $f \in \ker T_g$ with $f(0)=0$, then $gf \in 
\overline{H^2_0}$ and so $g (\overline z f)\in 
\overline{H^2_0}$ also, with $\overline z f \in H^2$, which means that $\overline z f \in \ker T_g$ too.
Indeed, a similar argument shows that we may divide out each inner factor while remaining in the kernel.

Thus Hitt proved the following result.

\begin{thm}\label{thm:hitt}
The nearly $S^*$-invariant  subspaces have the form
$M= uK$, with $u \in M$ of unit norm, $u(0)>0$,
and $u$   orthogonal to all elements of $M$ vanishing at the origin,  $K$ an $S^*$-invariant subspace,
and the operator of  multiplication by $u$ is  isometric from 
$K$ into $H^2$. 
\end{thm}

Note that $K$ may be $H^2$ itself, as for example $\theta H^2$ is
nearly $S^*$-invariant if $\theta$ is an inner function with $\theta(0) \ne 0$. This case is often overlooked, but these 
spaces $\theta H^2$ are not Toeplitz kernels, since they are not invariant
under dividing by $\theta$. The case we are most interested in is $K=K_\theta$, with
$\theta$   inner.

The link with $H^2(\AA)$ is that
if $M$ is an invariant subspace of $H^2(\AA)$, then under the change of variable $s=1/z$, the subspace $M \cap H^2(\CC \setminus \overline\DD)$ corresponds to a nearly $S^*$-invariant subspace. \\

Sarason gave a new proof of Hitt's theorem
using the de Branges--Rovnyak  spaces studied in  \cite{deBR}. He further showed that the inner function $\theta$ 
in the representation $\ker T_g= uK_\theta$ divides $(F-1)/(F+1)$, where $F$ is the Herglotz integral of $|u|^2$.

\subsection{The 1990s}

Hayashi
 \cite{hayashi90} and Sarason \cite{sarason94} continued to examine the
 nearly $S^*$-invariant subspaces   which are kernels of Toeplitz operators. 

Hayashi gave a complete characterization of such $uK_\theta$, as follows.
Let $u \in H^2$ be outer with $u(0)>0$, let $F$ be the Herglotz integral of $|u|^2$, and $b=(F-1)/(F+1)$.
Let $a$ be the outer function with $a(0)>0$ such that $|a|^2+|b|^2=1$ a.e.
We have $a=2f/(F+1)$ and $f=a/(1-b)$, and we write $u_\theta= a/(1-\theta b)$.

\begin{thm}
Let $M=uK_\theta$ as in Theorem \ref{thm:hitt}.
Then $M$ is the kernel of a Toeplitz operator if and only if $u$ is outer
and $a/(1-z\overline \theta b))^2$ is an exposed point of  the unit ball of $H^1$.
\end{thm}

Another way of writing this is to say that 

\begin{thm}
The nontrivial kernels of Toeplitz operators
are the subspaces of the the form
$M=u_\theta K_{z\theta}$, where $\theta$ is inner and $u \in H^2$ is outer with $u(0)>0$ and
$u^2$ an exposed point of the unit ball of $H^1$.
\end{thm}

Sarason gave an alternative proof of Hayashi's result, and a further
 discussion of rigid functions (for example the 1-dimensional Toeplitz kernels are
spanned by functions $u$ with $u^2$ rigid, and an outer function $u$ is rigid if and only 
if $\ker T_{\overline u/u}=\{0\}$) .

\subsection{The 2000s and 2010s}

Dyakonov \cite{Dyak00} took an alternative approach to Toeplitz kernels, using
 Bourgain's factorization for a unimodular function $\psi$ \cite{barclay,bourgain},
namely that
there is a triple $(B,b,g)$ such that $\psi=b\overline g/(Bg)$, where $b$ and $B$ are Blaschke products and $g$ is an invertible element in $H^\infty$.

As a result he showed the following result (in fact he showed a similar result in $H^p$ for $p > 1$).

\begin{thm}
 For every $\psi \in L^\infty\setminus\{0\}$, there exists a triple $(B,b,g)$ such that 
$\ker T_\psi=gb^{-1}(K_B \cap bH^2)$. 
\end{thm}

Then Makarov and Poltoratski \cite{MP05}, working in the upper half-plane $\CC^+$,
considered uniqueness sets. A Blaschke  set $\Lambda \subset \CC^+ $ is said to be a uniqueness set for $K_\theta$
 if every function in $K_\theta$ that vanishes on $\Lambda$ vanishes identically.  This property is  equivalent to the injectivity property for Toeplitz operators, i.e.,
$\ker T_{\overline\Theta B}=\{0\}$, where $B$ is the Blaschke product with zero set $\Lambda$.
Using these ideas they gave
a necessary and sufficient condition for the injectivity of a Toeplitz operator with the symbol $U=e^{i\gamma}$ where $\gamma$ is a real-analytic real function.

Before describing more recent work, we mention the survey article
of Hartmann and   Mitkovski \cite{HM16} and the book of Fricain and Mashreghi \cite{FricMash},
which give good treatments of the material we have discussed above.
Then the theory of model spaces and their operators (including composition operators, multipliers, restricted
shifts and indeed more general truncated Toeplitz operators) forms the subject of a monograph
\cite{GMR}. \\

\section{Near invariance and minimal kernels}

Toeplitz kernels form one of the most important classes of nearly $S^*$-invariant subspaces.
One may look at this property as meaning that if there is an element of a Toeplitz kernel $K$ of the
form $zf_+$ with $f_+ \in H^2$, then $f_+ \in K$.
In particular, one cannot have a one-dimensional Toeplitz kernel whose elements all vanish at $0$.
It is easy to see that an analogous property holds when $z$ is replaced by the inverse of a function $\eta \in \overline{H^\infty}$, as, for instance, an inner function.
More generally, if $\eta$ is a complex-valued function defined a.e. on $\TT$, we say that a proper closed subspace $\E$ of $H^2$ is {\em nearly $\eta$-invariant\/} if, for all $f_+ \in \E$, $\eta f_+ \in H^2$ implies that $\eta f_+ \in \E$. Thus,
saying that $\E$ is nearly $S^*$-invariant is equivalent to saying that $\E$ is nearly $\overline z$-invariant.

It can be shown \cite{CP14} that if $\eta \in H^\infty$ and $\eta$ is not constant, then
no finite-dimensional kernel is nearly $\eta$-invariant. However, one can characterise a vast class of functions
$\eta$, besides those in $\overline{H^\infty}$, for which all Toeplitz kernels are nearly
$\eta$-invariant. Let $\N_2$ denote the class of all such functions. We have the following.

\begin{thm}[\cite{CP14}]
If $\eta: X \to \CC$, measurable and defined on a set $X \subset \TT$ such that $\TT \setminus X$ has measure zero,
satisfies 
\[
L^2(\TT) \cap \eta \overline{H^2_0} \subset \overline{H^2_0},
\]
then every Toeplitz kernel is nearly $\eta$-invariant, i.e., $\eta \in \N_2$.
\end{thm}

Note that the class described in this theorem is rather large, including various well-known classes of functions,
not necessarily bounded \cite{CP14}, in particular all rational functions whose poles are in the closed disc $\overline\DD$
and all functions belonging to $\overline{H^2_0}$, as for instance those in $\overline\theta K_\theta= \overline z \overline{K_\theta}$
for some inner function $\theta$. We conclude therefore that if $\ker T_g \not= \{0\}$ (with $g \in L^\infty(\TT)$),
then, for each $\eta$
in that class, all $H^2$ functions that can be obtained from $f_+ \in \ker T_g$ by factoring out $\eta^{-1}$  must also belong to $\ker T_g$.
This establishes some sort of ``lower bound" for the Toeplitz kernel. For example, we have the following.

\begin{thm}[\cite{CP14}] A Toeplitz kernel that contains an element of the form $\phi_+=Rf_+$, where $f_+ \in H^2$ and $R \in H^\infty$ is a rational function of the form $R=p_1/p_2$, with $p_1$ and $p_2$ polynomials with no common zeroes, and $\deg p_1 \le \deg p_2$, has dimension at least $d:=P-Z+1$, where
$P$ is the number of poles of $R$, and $Z$ is the number of zeroes of $R$   in the exterior of  $\DD$.
\end{thm}

As another example, we have that if an inner function $\theta$ belongs to a Toeplitz kernel $K$, then $K \supset K_\theta$
\cite{CP14}. Thus, if $\theta$ is a singular inner function, then $K$ must be infinite-dimensional.

These lower bounds imply that, if $f_+ \in H^2$ has a non-constant inner factor, then
$\spam\{f_+\}$ cannot be a Toeplitz kernel. On the other hand, it is easy to see that there always exists a Toeplitz
kernel containing $f_+$, namely $\ker T_{\overline{f_+}/f_+}$, where the symbol is unimodular.
We are thus led to the question whether there is some ``smaller" Toeplitz kernel containing $f_+$. Or, in finite-dimensional language,
is there a minimum dimension for a Toeplitz kernel containing $f_+$? And can there be two different
Toeplitz kernels with that minimum dimension, such that $f_+$ is contained in both?
The answer to the first question is affirmative, while the second question has a negative answer. We have the
following result.

\begin{thm}[\cite{CP14}] \label{thm:3.3} Let $f_+ \in H^2 \setminus \{0\}$ and let $f_+=IO_+$ be its
inner--outer factorization. Then there exists a minimal Toeplitz kernel containing
$\spam\{f_+\}$, written $\Kmin(f_+)$, such that every Toeplitz kernel $K$ with $f_+ \in K$ contains 
$\Kmin(f_+)$, and we have
\beq\label{eq:3.1}
\Kmin(f_+)= \ker T_{\overline z \overline I \overline{O_+}/O_+}.
\eeq
\end{thm}

For example, given an inner function $\theta$, every kernel containing $\theta$ must contain $K_\theta$,
as mentioned before; the minimum kernel for $\theta$ is
\[
\Kmin(\theta)= \ker T_{\overline z\overline\theta}=K_\theta \oplus \spam\{\theta\}= K_\theta \oplus \theta K_z.
\]
If a Toeplitz kernel is the minimal kernel for $f_+ \in H^2$, we say that
$f_+$ is a {\em maximal function\/} or {\em maximal vector\/} for $K$. Since every Toeplitz kernel is the kernel of an operator
$T_{\overline z \overline I \overline{O_+}/O_+}$ for some inner function $I$ and outer
function $O_+ \in H^2$ \cite{sarason94} we conclude:

\begin{cor}
Every Toeplitz kernel has a maximal function.
\end{cor}

Note that this implies that every Toeplitz kernel $K$ contains an outer function, since, with the notation above, if $IO_+ \in K$,
then $O_+ \in K$ by near invariance.

One may ask when $\Kmin(f_+)=\spam\{f_+\}$, i.e., it is one-dimensional. There is a close connection between
one-dimensional Toeplitz kernels in $H^2$ and rigid functions in $H^2$.
It is easy to see that every rigid function is outer, and every rigid function in $H^1$ is the square of an outer function in $H^2$. We have the following.

\begin{thm}[\cite{sarason94}] If $f \in H^2 \setminus \{0\}$, then $\E=\spam\{f_+\}$ is a
Toeplitz kernel if and only if $f$ is outer and $f^2$ is rigid in $H^2$. In that case
$\E=\ker T_{\overline{f_+}/f_+}$.
\end{thm}

\section{Maximal functions in model spaces}

The maximal vectors for a given Toeplitz kernel can be characterized as follows.

\begin{thm}[\cite{CP17}] \label{thm:4.1} Let $g \in L^\infty \setminus \{0\}$ be such that $\ker T_g$ is nontrivial.
Then $k_+$ is a maximal vector for $\ker T_g$ if and only if
$k_+ \in H^2$ and $k_+=g^{-1}\overline z \overline{p_+}$, where $p_+ \in H^2$ is outer.
\end{thm}

Since model spaces are Toeplitz kernels ($K_\theta=\ker T_{\overline\theta}$), the maximal vectors are
the function $k_+ \in H^2$ of the form
\[
k_+ = \theta \overline z \overline{p_+} \qquad (p_+  \in H^2, \hbox{ outer}),
\]
i.e., such that $\theta \overline z \overline{k_+}$ is an outer function. Thus, the
reproducing kernel function, defined for each $w \in \DD$ by
\[
k_w^\theta(z):= \frac{1-\overline{\theta(w)}\theta(z)}{1-\overline w z}, \qquad (z \in \TT),
\]
is not in general a maximal vector for $K_\theta$, since
\[
\theta \overline z \overline{k^\theta_w} = \frac{\theta-\theta(w)}{z-w},
\]
which is not outer  in general. On the other hand, we have that
\[
\widetilde{ k_w^\theta}(z) :=  \frac{\theta(z)-\theta(w)}{z-w}  
\] is a maximal vector for $K_\theta$, for every $w \in \DD$.

Other maximal vectors for the model space $K_\theta$ can be found using the result that follows.
We use the notation $\G H^\infty$ for the set of invertible elements of the algebra $H^\infty$.

\begin{thm} \label{thm:4.2}
If $f_+$ is a maximal vector for $\ker T_g$, where $g \in L^\infty(\TT)$, then
$\theta h_+^{-1}f_+$ is a maximal vector for $\ker T_{h_-\overline\theta gh_+}$, for every
inner function $\theta$ and every $h_+ \in \G H^\infty$, $h_- \in \G\overline{H^\infty}$.
\end{thm}

\beginpf
From Theorem \ref{thm:4.1}, if $\Kmin(f_+)=\ker T_g$, then $gf_+=\overline z \overline{p_+}$,
where $p_+ \in H^2$ is outer. Therefore $\theta h^{-1}_+ f_+ \in H^2$ is such that
\[
h_- \overline\theta h_+ g (\theta h_+^{-1} f_+) = h_- gf_+ = \overline z(h_- \overline{p_+}),
\]
and using Theorem \ref{thm:4.1} again, we conclude that $\Kmin(\theta h^{-1}_+ f_+)=
\ker T_{h_- \overline\theta h_+ g}$.
\endpf

If the inner function is a finite Blaschke product $B$, with $B(z_0)=0$ for some $z_0 \in \DD$, then
it is easy to see from Theorem \ref{thm:3.3} that
\[
\Kmin\left( \frac{B}{z-z_0} \right) = \ker T_{\overline B}=K_B.
\]
Now each inner function $\theta$ can be factorized as
\[
\theta=h_- B h_+,
\]
where $B=\frac{\theta-a}{1-\overline a \theta}$ with $|a|<1$ is a Blaschke product and
$h_-=1+a\overline B \in \G \overline{H^\infty}$, and
$h_+ = \frac{1}{1+\overline a B} = \overline{h_-^{-1}} \in \G H^\infty$ \cite{nik};
thus it follows from Theorem \ref{thm:4.2} that
\beq\label{eq:4.6}
\phi^\theta_+ := \overline{h_-^{-1}} \frac{B}{z-z_0} = h_+ \frac{B}{z-z_0} = \frac{h_-^{-1}\theta}{z-z_0}
\eeq
is a maximal vector for $K_\theta= \ker T_{\overline\theta}$.

Note that, from \eqref{eq:4.6}, we can express $\theta$ in terms of these maximal vectors for $K_\theta$, using the same notation as above:
\beq\label{eq:4.7}
\theta=(z-z_0) h_- \phi^\theta_+.
\eeq

From Theorem \ref{thm:4.2}, applied to Toeplitz kernels that are model spaces, we also obtain the following.

\begin{thm} [\cite{CMP}]
Let $\theta$ and $\theta_1$ be inner functions. If $k_{1+}$ is a maximal vector for $K_{\theta_1}$,
then $\theta k_{1+}$ is a maximal vector for $K_{\theta\theta_1}= K_{\theta_1} \oplus \theta_1 K_\theta$.
\end{thm}

Thus if $\Kmin(k_{1+})$ is a model space $K_{\theta_1}$, then $\Kmin(\theta k_{1+})$
is also a model space, $K_{\theta\theta_1}$ for all inner functions $\theta$.

More generally, one can consider the minimal kernel containing a given set of functions.
In particular, when these functions are maximal vectors for model spaces, we obtain the
following generalization of the previous result.

\begin{thm} [\cite{CMP}]
Let $k_{1+}, k_{2+}, \ldots, k_{n_+} \in H^2$ be maximal vectors for $K_{\theta_1},K_{\theta_2},\ldots,K_{\theta_n}$,
respectively, where $\theta_j$ is an inner function for $j=1,2,\ldots,n$. Then there exists a minimal kernel
containing $\{k_{j+}: j=1,2,\ldots,n\}$, and for $\theta=\LCM(\theta_1,\theta_2,\ldots,\theta_n)$
we have
\[
K=K_\theta= \clos_{H^2} (K_{\theta_1}+K_{\theta_2}+ \cdots + K_{\theta_n})= K_{\theta_j }\oplus \theta_j K_{\theta \overline{\theta_j}},
\]
for each $j=1,2,\ldots,n$.
\end{thm}

\section{On the relations between $\ker T_g$ and $\ker T_{\theta g}$}
\label{sec:5}

Direct sum decompositions of the form $K_{\theta\theta_1}=K_{\theta_1}\oplus \theta_1 K_\theta$
can also be expressed in terms of maximal functions, using \eqref{eq:4.7} with $\theta$ replaced by $\theta_1$:

\beq\label{eq:5.1}
K_{\theta \theta_1} = K_{\theta_1} \oplus (z-z_0) h_- \phi_+^{\theta_1} K_\theta.
\eeq
For $g=\overline\theta \overline{\theta_1}$ the identity \eqref{eq:5.1} is equivalent to
\beq\label{eq:5.1A}
\ker T_g = \ker T_{\theta g} \oplus (z-z_0) h_- \phi_+^{\theta g}K_\theta,
\eeq
where $\phi_+^{\theta g}$ is a maximal vector for $\ker T_{\theta g}$ and $h_-=1$ if $\theta$ is
a Blaschke product with $\theta(z_0)=0$. This relation can be extended for general
$g \in L^\infty(\TT)$ when $\theta$ is a finite Blaschke product, in terms
of maximal functions and model spaces.

Indeed for every $g \in L^\infty(\TT)$ and every non-constant inner function $\theta$, we have
\[
\ker T_{\theta g} \subsetneq \ker T_g,
\]
whenever $\ker T_g \ne \{0\}$.

If $\theta$ is not a finite Blaschke product and $\dim \ker T_g < \infty$, then $\ker T_{\theta g}=\{0\}$; while,
if $\ker T_g$ is infinite-dimensional, then $\ker T_{\theta g}$ may or may not be finite-dimensional, and
in particular it can be $\{0\}$ -- as it happens, for instance, when
$\overline g$ is an inner function dividing $\theta$, or in the case of the following example.

\begin{ex} [\cite{CP17,CMP}] \label{ex:5.1}
{\rm
For $\theta(z)=\exp\left( \frac{z+1}{z-1} \right)$ and $\psi(z)=\exp\left( \frac{z-1}{z+1} \right)$, we have
$\ker T_{\overline z \theta \overline \psi}=\{0\}$.
}
\end{ex}

For finite Blaschke products $\theta$ we have the following.

\begin{thm} [\cite{CMP}] \label{thm:5.2} If $g \in L^\infty(\TT)$ and $\theta$ is a  finite Blaschke product, then
\[
\dim\ker T_g < \infty \quad \hbox{if and only if}\quad  \dim \ker T_{\theta g}< \infty,
\]
and $\ker T_g$ is finite-dimensional if and only if there
exists $k_0 \in \ZZ$ such that $\ker T_{z^{k_0}g}=\{0\}$; in that case $\dim \ker T_g \le \max \{0,k_0 \}$.
Moreover, if $\dim \ker T_g < \infty$, we have
\beq \label{eq:5.2}
\dim\ker T_{\theta g} = \max \{0, \dim\ker T_g - k \},
\eeq
where $k$ is the number of zeroes of $\theta$ counting their multiplicity.
\end{thm}

Thus, in particular, if $\dim\ker T_g = d < \infty$ and $\theta$ is 
a finite \BP such that $\dim K_\theta =k \le d$, then
\beq\label{eq:5.3}
\dim \ker T_{\theta g} = \dim \ker T_g - k.
\eeq
Of course, when $\ker T_g$ is infinite-dimensional and the same happens with $\ker T_{\theta g}$, it is not possible
to relate their dimension as in \eqref{eq:5.3}.
We can, however, use maximal functions to present an alternative relation, analogous to \eqref{eq:5.1A}, which not only generalizes Theorem
\ref{thm:5.2} but moreover sheds new light on the meaning of \eqref{eq:5.2} when $k < \dim \ker T_g < \infty$.

\begin{thm} [\cite{CMP}] Let $g \in L^\infty(\TT)$ and let $B$ be \ finite \BP of degree $k$. If $\dim\ker T_g \le k$, then $\ker T_{Bg}=\{0\}$; if $\dim\ker T_g > k$, then
\[
\ker T_g = \ker T_{Bg} \oplus (z-z_0) \phi_+ K_B,
\]
where $z_0$ is a zero of $B$ and $\phi_+$ is a maximal function for $\ker T_{Bg}$.
\end{thm}

\section{Injective Toeplitz operators}
\label{sec:6}

Clearly, the existence of maximal functions and the results of the previous section are closely connected with the question of injectivity of Toeplitz
operators, which in turn is equivalent to the question whether the Riemann--Hilbert problem $gf_+=f_-$, with $f_+ \in H^2$ and $f_- \in \overline{H^2_0}$, has a nontrivial solution.

It is well known that various properties of a Toeplitz operator, and in particular of its kernel, can be
described in terms of an appropriate factorization of its symbol
(\cite{BS, Dud, GK, LS, MP}).
For instance, the so-called $L^2$-factorization is a representation of 
the symbol $g \in L^\infty(\TT)$ as a product
\beq\label{eq:6.1}
g=g_- d g_+^{-1},
\eeq
where $g_+^{\pm 1} \in H^2$, $g_-^{\pm 1} \in \overline{H^2}$ and $d(z)=z^k$ for some $k \in \ZZ$.
If $g$ is invertible in $L^\infty(\TT)$ and admits an $L^2$-factorization, then 
$\dim\ker T_g=|k|$ if $k \le 0$, and $\dim \ker T_g^*=k$ if $k >0$.
The factorization \eqref{eq:6.1} is called a bounded factorization when
$g_+^{\pm 1},\overline{g_-^{\pm 1}} \in H^\infty$. In various subalgebras of $L^\infty(\TT)$, every
invertible element admits a factorization of the form \eqref{eq:6.1}, where the middle factor $d$ is an inner function. 
This is the case in the 
Wiener algebra on $\TT$ and in the analogous   
algebra $APW$  of almost-periodic functions
on the real line $\RR$.
In the latter case $d$ may be a singular inner function, $d(\xi)=\exp(-i\lambda \xi)$ with
$\lambda \in \RR$  and we have that if $g \in APW$ is invertible in $L^\infty(\RR)$ then $\ker T_g$ is either  trivial or isomorphic to
an infinite-dimensional model space $K_\theta$ with $\theta(\xi)=\exp(i \lambda \xi)$, depending on whether
$\lambda \le 0$ or $\lambda > 0$. For more details see \cite{CP17} and \cite[Sec. 8.3]{BKS}.

For $g_1,g_2 \in L^\infty(\TT)$, we say that $g_1 \sim g_2$ if and only if there are functions $h_+ \in \G H^\infty$, $h_- \in \G \overline{H^\infty}$
such that $g_1=h_- g_2 h_+$, and in that case we have $\ker T_{g_1} = h^{-1}_+ \ker T_{g_2}$
(which we write as $\ker T_{g_1} \sim \ker T_{g_2}$).
Thus if \eqref{eq:6.1} is a bounded factorization, we have $g \sim z^k$ and $\ker T_g= \{0\}$
if $k \ge 0$, and $\ker T_g \sim K_{z^{|k|}}$ if $k < 0$.

$L^2$ factorizations are  a particular case of factorizations of the form
\beq\label{eq:6.2}
g=g_- \theta^{-N} g_+^{-1}, \qquad g_- \in \overline{H^2}, \quad g_+ \in H^2,
\eeq
where $\theta$ is an inner function and $N \in \ZZ$. We have the following.

\begin{thm}[\cite{CP16,CP17}]
If $g \in L^\infty(\TT)$ admits a factorization \eqref{eq:6.2},where $\overline{g_-}$ and $g_+$ are outer functions in $H^2$,
with $g_+^2$ rigid in $H^1$, then
\[
\ker T_g \ne \{0\} \quad \hbox{if and only if} \quad N > 0.
\]
If $N>0$ and $\theta$ is a finite \BP of degree $k$, then $\dim\ker T_g=kN$; if $\theta$ is not a finite 
Blaschke product, then $\dim \ker T_g=\infty$.
\end{thm}

We also have the following.

\begin{thm} [\cite{CP16,nak}]
For $g \in L^\infty(\TT)$, $\ker T_g$ is nontrivial of finite dimension if and only if, for some $N \in \NN$, $g$
admits a factorization $g=g_- z^{-N} g_+^{-1}$, where $g_- \in \overline{H^2_0}$ is outer, 
and $g_+ \in H^2$ is outer with $g^2_+$ rigid in $H^1$.
In that case $\ker T_g = \ker T_{z^{-N}\overline{g_+}/g_+}$, and $\dim\ker T_g=N$.
\end{thm}
Some other results regarding conditions for injectivity or non-injectivity of
Toeplitz operators will be mentioned in the next section.

\section{Multipliers between Toeplitz kernels}

The existence of maximal vectors for every non-zero Toeplitz kernel also provides test functions for
various properties of these spaces. 

In \cite{crofoot} Crofoot characterized the multipliers
from a
model space {\em onto\/} another. Partly motivated by that work,
Fricain, Hartmann and Ross addressed in 
\cite{FHR}
the question of which holomorphic functions $w$ multiply a model space $K_\theta$ into another model space $K_\phi$.
Their main result shows that
$w$ multiplies $K_\theta$ into $K_\phi$
(written $w \in \M(K_\theta,K_\phi)$) if and only if
\begin{enumerate}[(i)]
\item $w$ multiplies the function $S^*\theta= \widetilde k_0^\theta$ into $K_\phi$, and
\item $w$ multiplies $K_\theta$ into $H^2$, which can be expressed by saying that
$|w|^2 \, dm$ is a Carleson measure for $K_\theta$.
\end{enumerate}

Model spaces being a particular type of Toeplitz kernel, that question may be posed more generally for the latter. We may
also ask whether more general test functions can be used, other than $S^* \theta$.

In this more general setting, one immediately notices that, unlike multipliers between model spaces,
multipliers between general Toeplitz kernels need not lie in $H^2$.
In fact, for model spaces, we must have $w \in H^2$ if $w \in \M(K_\theta,K_\phi)$, because we must then have
$w k^\theta_0 \in K_\phi \subset H^2$, and $1/k^\theta_0  \in H^\infty$; but the
function $w(z)=(z-1)^{-1/2}$ multiplies $\ker T_g$, with
$g(z)=z^{-3/2}$ and $\arg z \in [0,2\pi)$ for $z \in \TT$, onto the model
space $K_z=\ker T_{\overline z}$ consisting of the constant functions, even though $w \not\in H^2$.

One can characterize all multipliers from one Toeplitz kernel into another as follows. We denote
by $\C(\ker T_g)$ the class of all $w$ such that $|w|^2 \, dm$ is a Carleson measure for $\ker T_g$, i.e., $w \ker T_g \subset L^2(\TT)$, and by $\N_+$ the Smirnov class.

\begin{thm} [\cite{CP17}] \label{thm:7.1}
Let $g,h \in L^\infty(\TT) \setminus \{0\}$ be such that $\ker T_g$ and $\ker T_h$ are nontrivial.
Then the following are equivalent:
\begin{enumerate}[(i)]
\item $w \in \M(\ker T_g, \ker T_h)$;
\item $w \in \C(\ker T_g)$ and $wk_+ \in \ker T_h$ for some (and hence all) maximal vectors $k_+$
of $\ker T_g$;
\item $w \in \C(\ker T_g)$ and $hg^{-1}w \in \overline{\N_+}$.
\end{enumerate}
\end{thm}

Note that if $k_+$ is not a maximal vector for $\ker T_g$, then $k_+$ cannot be used
as a test function; for example, the function $w=1$ is not a multiplier from $\ker T_g$
into $\Kmin(k_+)$, even though $wk_+ \in \Kmin(k_+)$.

\begin{cor}[\cite{CP17}]
With the same assumptions as in Theorem \ref{thm:7.1}, and assuming moreover that $hg^{-1} \in L^\infty(\TT)$, one has
\[
w \in \M(\ker T_g, \ker T_h) \quad \hbox{if and only if} \quad w \in \C(\ker T_g) \cap \ker T_{\overline z gh^{-1}}.
\]
\end{cor}

By considering the special case $g = \overline \theta$, where $\theta$ is inner, we obtain the
following result.

\begin{cor}[\cite{CP17,FHR}]
Let $\theta$ be inner and let $h \in L^\infty(\TT) \setminus \{0\}$ be such that
$\ker T_h$ is nontrivial. Then the following are equivalent:
\begin{enumerate}[(i)]
\item $w \in \M(K_\theta, \ker T_h)$;
\item $w \in \C(K_\theta)$ and $wS^*\theta \in \ker T_h$;
\item $w \in \C(K_\theta) \cap \ker T_{\overline z\theta h}$.
\end{enumerate}
\end{cor}

The last two corollaries also bring out a close connection between
the existence of non-zero multipliers in $L^2(\TT)$ and their description, on the one hand,
and the question of injectivity of an associated Toeplitz operator and the characterization of its kernel
(discussed in Sections \ref{sec:5} and \ref{sec:6}), on the other hand.
Thus, for instance, the result of Example \ref{ex:5.1} implies that, since
$T_{\overline z \theta \overline \psi}$ is injective in that case, we have $\M(K_\theta,K_\phi)=\{0\}$.
Another example is the following:

\begin{ex}{\rm 
Let $\theta$, $\phi$ be two inner functions with $\phi \preceq \theta$, i.e., $K_\phi \subset K_\theta$.
Then $\dim \ker T_{\overline z \theta \overline \phi} \le 1$, since $\theta\overline\phi \in H^\infty$ and
$\ker T_{\theta\overline \phi} = \{0\}$
(\cite{BCD}).
We have $\ker T_{\overline z\theta \overline \phi}=\CC$ if $\phi=a\theta$ with $a \in \CC$, $|a|=1$,
and we have $\ker T_{\overline z\theta \overline \phi}=\{0\}$ if $\phi \prec \theta$; therefore
$\M(K_\theta,K_\phi) \ne \{0\}$
if and only if $K_\theta=K_\phi$, in which case $\M(K_\theta,K_\phi)=\CC$.
}
\end{ex}

The class of bounded multipliers,
\[
\M_\infty(\ker T_g, \ker T_h)=\M(\ker T_g , \ker T_h) \cap H^\infty,
\]
 is of great importance.
For instance, the question whether $w=1$ is a multiplier from $\ker T_g$ into $\ker T_h$ is equivalent to asking whether $\ker T_g \subset \ker T_h$.
Noting that the Carleson measure condition is redundant for bounded
$w$, we obtain the following characterization from Theorem \ref{thm:7.1}.

\begin{thm}  \label{thm:7.5}
Let $g,h \in L^\infty(\TT) \setminus \{0\}$ be such that $\ker T_g$ and $\ker T_h$ are nontrivial.
Then the following are equivalent:
\begin{enumerate}[(i)]
\item $w\in \M_\infty (\ker T_g, \ker T_h)$;
\item $w \in H^\infty$ and $wk_+ \in \ker T_h$ for some (and hence all) maximal vectors $k_+$
of $\ker T_g$;
\item $w \in H^\infty$ and $hg^{-1}w \in \overline{H^\infty}$ (assuming that $hg^{-1} \in L^\infty(\TT)$).
\end{enumerate}
\end{thm}

For model spaces, we thus recover the main theorem on bounded multipliers from \cite{FHR}:

\begin{cor}
Let $\theta$ and $\phi$ be inner functions, and let $w \in H^2$. Then
\begin{eqnarray*}
w \in \M_\infty(K_\theta,K_\phi) &\iff& w \in \ker T_{\overline z \theta \overline \phi} \cap H^\infty
\iff wS^*\theta \in K_\phi \cap H^\infty \\
&\iff& w \in H^\infty \quad \hbox{and} \quad  \theta\overline\phi w \in \overline{H^\infty}.
\end{eqnarray*}
\end{cor}

Applying the results of Theorem \ref{thm:7.5} to $w=1$ we obtain moreover the following results.

\begin{cor}
Under the same assumptions as in Theorem \ref{thm:7.5}, the following conditions are equivalent:
\begin{enumerate}[(i)]
\item $ \ker T_g \subset \ker T_h$;
\item $hg^{-1} \in \overline{\N_+}$;
\item there exists a maximal function $k_+$ for $\ker T_g$ such that $k_+ \in \ker T_h$.

{\hskip-25pt If, moreover, $\ker T_g$ contains a maximal vector $k_+$ with $k_+, k_+^{-1} \in L^\infty(\TT)$,}

{\hskip-25pt then each of the above conditions is equivalent to}
\item $k_+ \in \ker T_h \cap H^\infty$.
\end{enumerate}
\end{cor}

\begin{cor}
Under the same assumptions as in Theorem \ref{thm:7.5},
if $hg^{-1} \in \G L^\infty(\TT)$, then
\[
\ker T_g \subset \ker T_h \quad \hbox{if and only if} \quad  hg^{-1} \in \overline{H^\infty}.
\]
\end{cor}

This last result implies in particular that, assuming that $hg^{-1} \in L^\infty(\TT)$,  
a Toeplitz kernel is contained in another Toeplitz kernel
if and only if they take the form $\ker T_g$ and $\ker T_{\theta g}$ for some inner function
$\theta$ and $g \in L^\infty(\TT)$ (cf. Section \ref{sec:5}).

\begin{cor}\label{cor:7.9}
Under the same assumptions as in Theorem \ref{thm:7.5},
we have $\ker T_g=\ker T_h$ if and only if $\displaystyle \frac{g}{h} = \frac{\overline{p_+}}{\overline{q_+}}$
with $p_+, q_+ \in H^2$  outer.
If moreover $hg^{-1} \in \G L^\infty(\TT)$, then we have
\[
\ker T_g = \ker T_h \quad \hbox{if and only if} \quad hg^{-1} \in \G \overline{H^\infty}.
\]
\end{cor}

We can draw several interesting conclusions from these results.

1. First, we can characterize the Toeplitz kernels that are contained in a given model space $K_\theta$
($\ker T_g = \ker T_{\overline\theta \alpha}$, with $\alpha$ inner), and those that contain $K_\theta$
($\ker T_g$ with $g \in \overline{\theta H^\infty}$), assuming that the symbols are in $\G L^\infty(\TT)$.

2. Second, while \eqref{eq:3.1} provides an expression for a (unimodular) symbol $g$ such that $\ker T_g$ is the minimal kernel
for a given function with inner--outer factorization
$\phi_+= IO_+$, it is not claimed that all Toeplitz operators with that
kernel have the same symbol. Indeed, from Corollary \ref{cor:7.9}, we have that if
$\ker T_g= \Kmin(\phi_+)$ with $\phi_+=IO_+$, then
$g=\displaystyle \frac{\overline{p_+}}{\overline{q_+}}  \frac{\overline I\overline{O_+}}{O_+}$
with 
$p_+, q_+ \in H^2$  outer; if, moreover, $g \in \G L^\infty(\TT)$, then
$g= h_- \overline I \overline{O_+}/\overline O_+$, with $h_- \in \G \overline{H^\infty}$.

3. It is clear that a Toeplitz operator with unimodular symbol $u$ is non-injective if and only if it has a maximal
vector, i.e., there exist an inner function $I$ and an outer function $O_+ \in H^2$
such that $\ker T_u = \Kmin(IO_+)= \ker T_{\overline I\overline{O_+}/O_+}$, which is equivalent,
as shown in point 2, to having
\[
u= \overline{z} \frac{\overline I\overline{O_+}}{O_+} \, h_-, \qquad \hbox{with} \quad h_- \in \G 
\overline{H^\infty}.
\]
Since $|h_-|=1$ a.e. on $\TT$, we conclude that $h_-$ must be a unimodular constant,
and therefore $T_u$ is non-injective if and only if
\[
u= \overline z \overline I \overline{O_+} / O_+,
\]
thus recovering a result by Makarov and Poltoratski
\cite[Lem.~3.2]{MP05}.

4. Since there are different maximal functions for each Toeplitz kernel
with dimension greater than 1, one may ask how they can be related.
Again, from Corollary \ref{cor:7.9}, we see that if
$\Kmin(f_{1+})=\Kmin(f_{2+})$, where $f_{1+}=I_1 O_{1+}$ and $f_{2+}=I_2 O_{2+}$
with $I_1,I_2$ inner and $O_{1+}, O_{2+} \in H^2$ outer, then
\[
\frac{\overline{I_1}\overline{O_{1+}}}{O_{1+}}
= \frac{\overline{I_2}\overline{O_{2+}}}{O_{2+}} \, h_-,
\]
where $h_- \in \G \overline{H^\infty}$, $|h_-|=1$, and so $h_-$ is constant. Thus, finally
$f_{1+}$ and $f_{2+}$ are related by
\[
f_{2+}=f_{1+} \frac{\overline{O_{2+}}}{\overline{O_{1+}}}.
\]

\subsection*{Acknowledgements}

This work was partially supported by FCT/Portugal through the grant \\ UID/MAT/04459/2013.

\end{document}